\documentclass[12pt,reqno]{amsart}
\usepackage{latexsym,amssymb,amsmath,multicol, rotating,lscape}
 2

\newtheorem{thm}{Theorem}[section]

\newtheorem{prop}[thm]{Proposition}

\newcommand{\thmref}[1]{Theorem~\ref{#1}}

\theoremstyle{remark}

\begin{document}

\title[Representation of certain quadratic forms]
{Certain quaternary quadratic forms of level 48 and their representation numbers}
 
\author{B. Ramakrishnan, Brundaban Sahu and Anup Kumar Singh}
\address[B. Ramakrishnan and Anup Kumar Singh]{Harish-Chandra Research Institute, HBNI,  
       Chhatnag Road, Jhunsi,
     Allahabad -     211 019,
   India.}
\address[Brundaban Sahu]
{School of Mathematical Sciences, National Institute of Science 
Education and Research, Bhubaneswar, HBNI, 
Via- Jatni, Khurda, Odisha 752 050  
India.}

\email[B. Ramakrishnan]{ramki@hri.res.in}
\email[Brundaban Sahu]{brundaban.sahu@niser.ac.in}
\email[Anup Kumar Singh]{anupsingh@hri.res.in}


\date{\today}

\begin{abstract}

In this paper, we find a basis for the space of modular forms of weight $2$ on $\Gamma_1(48)$. We use this basis to find formulas for the number of representations of a positive integer $n$ by certain quaternary quadratic forms of the form $\sum_{i=1}^4 a_i x_i^2$, $\sum_{i=1}^2 b_i(x_{2i-1}^2 + x_{2i-1}x_{2i}+x_{2i}^2)$ and 
$a_1x_1^2 + a_2 x_2^2 + b_1(x_3^2+x_3x_4+x_4^2)$, where $a_i$'s belong to $\{1,2,3,4,6,12\}$ and $b_i$'s belong to $\{1,2,4,8,16\}$.
\end{abstract}

\maketitle

\section{Introduction}
In this paper we consider the problem of finding the number of representations of a natural number by the following three types of quaternary quadratic forms given by  
\begin{equation*}
\begin{split}
{\mathcal Q}_1: & ~a_1x_1^2+a_2x_2^2+a_3x_3^2+a_4x_4^2,\\
{\mathcal Q}_2: & ~b_1(x_1^2+x_1x_2+x_2^2) +  b_2(x_3^2+x_3x_4+x_4^2),\\
{\mathcal Q}_3:& ~a_1x_1^2 + a_2 x_2^2 + b_1(x_3^2+x_3x_4+x_4^2),
\end{split}
\end{equation*}
where the coefficients $a_i\in \{1, 2, 3, 4, 6, 12\},~ 1 \leq i \leq 4$ and $b_i\in \{1,2,4,8,16\}$, $i=1,2$.  Without loss of generality we may assume that $a_1 \leq a_2 \leq 
a_3 \leq a_4$, $b_1\le b_2$ and $\gcd(a_1,a_2,a_3,a_4)=1$, $\gcd(b_1,b_2)=1$, $\gcd(a_1,a_2,b_1)=1$.  
Finding explicit formula for the number of representations of $n$ by these types of  quadratic forms is a classical problem in number theory. The classical formula of 
Jacobi for the sum of four squares correspond to the quadratic form ${\mathcal Q}_1$ with $(a_1,a_2,a_3,a_4) = (1,1,1,1)$. There are several works in the literature which consider giving formulas for the representation numbers corresponding to quaternary quadratic forms with coefficients. We list some of them here \cite{{aaa},{aalw},
{aalw1}, {aalw2}, {aalw3}, {aaw}, {aaw1}, {apw}, {aw}, {w}}. 
In \cite{aaa}, A. Alaca et. al considered 35 quadratic forms of type ${\mathcal Q}_1$ with coefficients $a_i\in \{1, 2, 3, 6\}$. 
They constructed explicit bases for the spaces of modular forms of weight $2$ on $\Gamma_0(24)$ with character $\chi_0$ (trivial character modulo $24$) or $\chi_d = \left(\frac{d}{\cdot}\right)$ for $d = 8, 12, 24$ and used these bases to determine formulas for the number of representations of a natural number by ${\mathcal Q}_1$, 
with $a_i\in \{1,2,3,6\}$.  However, out of these $35$ quadratic forms of type ${\mathcal Q}_1$, formulas for $15$ forms appeared in the papers 
\cite{{aalw2}, {aalw3}, {aw}}. Precisely, the cases $(1,1,1,4), (1,1,2,2), (1,1,3,3)$ are considered in \cite{aalw}, the cases $(1,1,1,3), (1,1,2,6), (1,2,2,3), (1,3,3,3), (1,3,6,6), (2,3,3,6)$ are considered in \cite{aalw2}, the cases $(1,1,1,2), (1,2,2,2)$ are considered in \cite{aalw3} and the cases $(1,1,2,3), (1,2,2,6), (1,3,3,6), (2,3,6,6)$ are considered in \cite{aw}. 

The total number of quadratic forms ${\mathcal Q}_1$, with $a_i\in \{1,2,3,4,6,12\}$ is $126$. Out of this, $36$ cases come from the works of \cite{{aaa}, {aalw}, 
{aalw2},{aalw3}, {aw}} (when $a_i\not=4,12$) and $35$ cases have the $\gcd$(coefficients) $>1$ and so in the present work we consider the remaining $55$ cases of this type of quadratic forms. There are only 4 quadratic forms  of type ${\mathcal Q}_2$ and there are $65$  quadratic forms of type ${\mathcal Q}_3$ (such that the 
coefficients  have no common factors). In the following table we give the list of quadratic forms ${\mathcal Q}_i$, $i=1,2,3$ considered in this paper. They are listed according to the modular forms space (in which the corresponding theta series belong). Note that in place of $M_2(48, \chi)$, we mention only the character $\chi$ which is either $\chi_0$ (trivial character modulo $48$) or $\chi_d$, $d=8,12,24$.\\
\begin{center}
Table 1. (List of quadratic forms)\\
\smallskip

\begin{tabular}{|c|c|} 
\hline
&\\
Space &  $ (a_1, a_2, a_3, a_4)$ for ${\mathcal Q}_1$ \\
&\\
\hline
&$(1, 1, 1, 4), (1, 1, 4, 4), (1, 1, 3, 12)$, $(1, 1, 12, 12), (1, 2, 2, 4),  (1, 2, 6, 12), $\\
$\chi_0$& $(1, 3, 3, 4), (1, 3, 4, 12),$ $(1, 4, 4, 4), (1, 4, 6, 6), (1, 4, 12, 12),$ \\
& $ (2, 2, 3, 12), (2, 3, 4, 6), (3, 3, 4, 4), (3, 4, 4, 12)$ \\
\hline
& $(1, 1, 2, 4), (1, 1, 6, 12), (1, 2, 4, 4), (1, 2, 3, 12), (1, 2, 12, 12), $  \\
$\chi_8$ & $(1, 3, 4, 6),  (1, 4, 6, 12), (2, 3, 3, 4),$ $ (2, 3, 4, 12), (3, 4, 4, 6)$\\
\hline
& $(1, 1, 1, 12), (1, 1, 3, 4), (1, 1, 4, 12), (1, 2, 2, 12), (1, 2, 4, 6), (1, 3, 3, 12),$\\
 $\chi_{12}$ & $(1, 3, 4, 4), (1, 3, 12, 12), (1, 4, 4, 12), (1, 6, 6, 12), (1, 12, 12, 12), (2, 2, 3, 4), $ \\
& $(2, 3, 6, 12), (3, 3, 3, 4), (3, 3, 4, 12), (3, 4, 4, 4), (3, 4, 6, 6), (3, 4, 12, 12)$ \\
 \hline
& $(1, 1, 2, 12), (1, 1, 4, 6), (1, 2, 3, 4), (1, 2, 4, 12), (1, 3, 6, 12), (1, 4, 4, 6),$ \\
 $\chi_{24}$ &$ (1, 6, 12, 12), (2, 3, 3, 12),  (2, 3, 4, 4), (2, 3, 12, 12), (3, 3, 4, 6), (3, 4, 6, 12)$ \\
\hline
\hline 
&\\
 &  $(b_1, b_2)$ for ${\mathcal Q}_2$ \\
&\\
\hline
$\chi_0$ & $(1,2), (1,4), (1,8), (1,16)$ \\
\hline 
\hline 
 &\\
 &  $(a_1,a_2, b_1)$ for ${\mathcal Q}_3$ \\
&\\
\hline
&$(1, 3, 1), (1,3, 2), (1, 3, 4), (1,3,8), (1, 3, 16), (1, 12, 1), (1, 12, 2), (1, 12, 4), (1, 12, 8),$\\  
$\chi_0$ & $(1, 12, 16), (2,6,1),(3,4,1),(3,4,2), (3,4,4), (3,4,8), (3,4,16), (4,12,1)$ \\ 
\hline
&  $(1,6,1),(1,6,2),(1,6,4),(1,6,8),(1,6,16), (2,3,1),$\\
$\chi_8$ & $(2,3,2),(2,3,4),(2,3,8),(2,3,16),(2,12,1),(4,6,1)$\\
\hline
& $(1,1,1),(1,1,2),(1,1,4),(1,1,8),(1,1,16),(1,4,1),(1,4,2),(1,4,4),$\\
$\chi_{12}$& $(1,4,8),(1,4,16), (2,2,1),(3,3,1),(3,3,2),(3,3,4),(3,3,8), (3,3,16),$\\ 
 & $(3,12,1), (3,12,2), (3,12,4), (3,12,8), (3,12,16),(4,4,1),(6,6,1),(12,12,1) $\\ 
\hline 
& $(1,2,1),(1,2,2),(1,2,4),(1,2,8),(1,2,16), (2,4,1), (3,6,1)$,\\
$\chi_{24}$ & $(3,6,2),(3,6,4),(3,6,8),(3,6,16), (6,12,1)$\\
\hline
\end{tabular}
\end{center}

\smallskip

Some of our formulas were also proved in the earlier works of K. S. Williams and his co-authors, which we  mention in the table below. However, our approach is different from their method. 

\bigskip

\newpage

\begin{center}
Table A. (List of earlier results)\\
\smallskip

\begin{tabular}{|l|c|c|} 
\hline
Type & Cases& Ref.\\
\hline 
 &$(1,1,1,12), (1,1,3,4), (1,1,4,12), (1,2,2,12), (1,3,3,12), (1,3,4,4)$ &\\
& $(1,3,12,12), (1,4,4,12), (1,6,6,12), (1,12,12,12), (2,2,3,4), (3,3,3,4)$, &  \cite{aalw1} \\
${\mathcal Q}_1$& $(3,3,4,12), (3,4,4,4), (3,4,6,6), (3,4,12,12)$, &\\
 & $(1,2,4,6)$ & \cite{aalw2}\\
 & $(1,1,2,4), (1,2,4,4)$ & \cite{aalw3}\\
\hline
${\mathcal Q}_2$ & $(1,2), (1,4)$ & \cite{aaw}\\
\hline
 & $(1,1,1), (1,1,2),  (1,1,4), (3,3,1), (3,3,2), (3,3,4)$ & \cite{aalw2}\\
${\mathcal Q}_3$ & $(1,1,8), (1,4,2), (3,12,2)$ & \cite{aaw1}\\
&$(1,3,1), (1,3,2),  (1,3,4), (1,4,4), (2,2,1), (6,6,1)$ & \cite{apw}\\
\hline
\end{tabular}
\end{center}

\bigskip

\smallskip

Let $\mathbb N, \mathbb N_0, \mathbb Z, \mathbb Q$ and $\mathbb C$ denote the sets of natural numbers, non-negative integers, integers, rational numbers and complex numbers respectively.  For $n \in \mathbb N,$ let us denote by
\begin{equation}
\begin{split}
{\mathcal N}_1(a_1, a_2, a_3, a_4; n)  &=\#\{ (x_1, x_2, x_3, x_4) \in \mathbb Z^4: a_1x_1^2+a_2x_2^2+a_3x_3^2+a_4x_4^2 = n \},\\
{\mathcal N}_2(b_1, b_2; n) ~=~ &\#\{ (x_1, x_2, x_3, x_4) \in \mathbb Z^4: b_1(x_1^2+x_1x_2+x_2^2) +  b_2(x_3^2+x_3x_4+x_4^2)  = n \},\\
{\mathcal N}_3(a_1, a_2, b_1; n) &= \#\{ (x_1, x_2, x_3, x_4) \in \mathbb Z^4: a_1x_1^2 + a_2 x_2^2 + b_1(x_3^2+x_3x_4+x_4^2) = n \}
\end{split}
\end{equation}
be the number of representations of $n$ by the quadratic forms ${\mathcal Q}_i$, $i=1,2,3$ respectively. In this paper, we construct explicit bases for the spaces 
of modular forms of weight $2$ on $\Gamma_0(48)$ with character $\chi$ (modulo $48$) and use them to give formulas for ${\mathcal N}_1(a_1,a_2,a_3,a_4;n)$,
${\mathcal N}_2(b_1,b_2;n)$ and ${\mathcal N}_3(a_1,a_2,b_1;n)$. It is to be noted that in his work \cite{ramanujan}, Ramanujan gave a list of $55$ universal quadratic forms of type ${\mathcal Q}_1$. Our work includes $8$ out of these $55$ forms which are given by $(a_1,a_2,a_3,a_4) = (1,1,1,4), (1,1,2,4), (1,1,2,12), (1,1,3,4),$  
$(1,2,3,4), (1,2,4,4), (1,2,4,6)$ and $(1,2,4,12)$. 

\section{Preliminaries and Statement of the Results}

We use the theory of modular forms to prove our results and so we first fix our notations and present some of the basic facts on modular forms. 
For positive  integers $k, N \ge 1$ and a Dirichlet character $\chi$ modulo $N$ with $\chi(-1) = (-1)^k$, let  
$M_k(N, \chi)$ denote the ${\mathbb C}$- vector space of holomorphic modular forms of weight $k$ for the congruence subgroup $\Gamma_0(N)$, with character $\chi$. Let us denote by $S_k(N, \chi)$,  the subspace of cusp forms in $M_k(N, \chi)$. When $\chi$ is the principal character modulo $N$, then we 
drop the symbol $\chi$ in the notation and write only $M_k(N)$ or $S_k(N)$.  The modular forms space is decomposed into the space of 
Eisenstein series (denoted by ${\mathcal E}_k(N, \chi)$) and the space of cusp forms $S_k(N, \chi)$ and one has 
\begin{equation}
M_k(N, \chi) =  {\mathcal E}_k(N, \chi) \oplus S_k(N, \chi).
\end{equation}
Explicit basis for the space ${\mathcal E}_k(N, \chi)$ can be obtained using the following construction. For details we refer to \cite{{miyake}, {stein}}. 
Suppose that $\chi$ and $\psi$ are primitive Dirichlet characters with conductors $N$ and $M$, respectively. For a positive integer $k \ge 2$, let 
\begin{equation}\label{eisenstein}
E_{k,\chi,\psi}(z) :=  c_0 + \sum_{n\ge 1}\left(\sum_{d\vert n} \psi(d)  \chi(n/d) d^{k-1}\right) q^n,
\end{equation}
where $q=e^{2 i\pi z} ~(z\in {\mathbb C}, {\rm Im}(z)>0)$ and 
$$
c_0 = \begin{cases}
0 &{\rm ~if~} N>1,\\
- \frac{B_{k,\psi}}{2k} & {\rm ~if~} N=1,
\end{cases}
$$
with $B_{k,\psi}$ denoting the generalized Bernoulli number with respect to the character $\psi$. 
Then, the Eisenstein series $E_{k,\chi,\psi}(z)$ belongs to the space $M_k(NM, \chi/\psi)$, provided 
\linebreak 
$\chi(-1)\psi(-1) = (-1)^k$ and $NM\not=1$. We give a notation to the inner sum in \eqref{eisenstein}:
\begin{equation}\label{divisor}
\sigma_{k-1, \chi,\psi}(n) := \sum_{d\vert n} \psi(d) \chi(n/d) d^{k-1}.
\end{equation}
When $\chi=\psi =1$ (i.e., when $N=M=1$) and $k\ge 4$, we have $E_{k,\chi,\psi}(z) = - \frac{B_k}{2k} E_k(z)$, where $E_k(z)$ is the normalized Eisenstein series of weight $k$ in the space $M_k(1)$, defined  by 
\begin{equation}
E_k(z) = 1 - \frac{2k}{B_k}\sum_{n\ge 1} \sigma_{k-1}(n) q^n. 
\end{equation}
In the above $\sigma_r(n)$ is the sum of the $r$-th powers of the positive divisors of $n$ and $B_k$ is the $k$-th Bernoulli number defined by $\displaystyle{\frac{x}{e^x-1} = \sum_{m=0}^\infty \frac{B_m}{m!} x^m}$.
We also need the Eisenstein series of weight $2$, which is a  quasimodular form of weight $2$, depth $1$ on $SL_2({\mathbb Z})$ and is given by 
\begin{equation}
E_2(z) = 1 -24\sum_{n\ge 1} \sigma(n) q^n. 
\end{equation}
(Note that $\sigma(n) = \sigma_1(n)$.) 
Let 
\begin{equation}
\eta(z)=q^{1/24} \prod_{n\ge1}(1-q^n)
\end{equation}
denote the Dedekind eta function. Then,  an eta-quotient is a finite product of integer powers of $\eta(z)$ and we denote it as follows:
\begin{equation}
\prod_{i=1}^s \eta^{r_i}(d_iz):=d_1^{r_1}d_2^{r_2} \cdots d_s^{r_s},
\end{equation}
where $d_i$'s are positive integers and  $r_i$'s are non-zero integers.
\smallskip

In the case of the space of cusp forms $S_k(N,\chi)$ we use a basis consisting of newforms of level $N$ and oldforms generated by the newforms of lower level $d$, $d\vert N$, $\chi$ modulo $d$, $d\not =N$. However, when $\chi =\chi_{12}$, we construct a basis for the space of newforms, which are not Hecke eigenforms. 
For a basic theory of newforms we refer to \cite{{a-l},{li}} and for details on modular forms, we refer to \cite{{koblitz},  {miyake}, {stein}}. 

\smallskip
\bigskip


We now state the main results of this paper. In the following statements $\chi$ denotes a Dirichlet character modulo $48$, which is either the principal character modulo 
$48$, denoted as $\chi_0$ or the Kronecker symbol $\chi_d = \left(\frac{d}{\cdot}\right)$, where $d=8,12,$ or $24$. For each such $\chi$, let $\ell_\chi$ denote the dimension of the ${\mathbb C}$- vector space $M_2(\Gamma_0(48),\chi)$. Then 
$$
\ell_\chi = \begin{cases}
14 & \mbox{if~} \chi = \chi_0 \mbox{~or~} \chi_{12},\\
12 & \mbox{if~} \chi = \chi_8 \mbox{~or~} \chi_{24}.\\
\end{cases}
$$

\begin{thm}\label{Q1}
Let $n\in \mathbb N$. For each entry $(a_1, a_2, a_3, a_4)$ corresponding to ${\mathcal Q}_1$ in Table 1, the associated theta series is a modular form of 
weight $2$ on $\Gamma_0(48)$ with character $\chi$. Therefore, using the basis given Table B (in \S 3.5), we have  
\begin{equation}
{\mathcal N}_1(a_1, a_2, a_3, a_4; n) = \sum_{i=1}^{\ell_\chi} \alpha_{i,\chi} A_{i,\chi}(n),
\end{equation}
where $A_{i,\chi}(n)$ are the Fourier coefficients of the basis elements $f_{i,\chi}$ and the values of the constants $\alpha_{i,\chi}$'s are given in (\S 6, Table 2). 
\end{thm}

\begin{thm}\label{Q2}
Let $n\in \mathbb N$.  Then we have 
\begin{equation}
\begin{split}
{\mathcal N}_2(1,2; n) & = 6 \sigma(n) - 12 \sigma(n/2) + 18 \sigma(n/3) - 36 \sigma(n/6),\\
{\mathcal N}_2(1,4; n) & = 6 \sigma(n) - 18 \sigma(n/2) - 18 \sigma(n/3) + 24 \sigma(n/4) + 54 \sigma(n/6) - 72 \sigma(n/12), \\
{\mathcal N}_2(1,8; n) & = \frac{3}{2} \sigma(n) - \frac{9}{2} \sigma(n/2) + \frac{9}{2} \sigma(n/3) + 9 \sigma(n/4) - \frac{27}{2} \sigma(n/6) - 12 \sigma(n/8)\\
& ~\quad  + 27 \sigma(n/12) - 36 \sigma(n/24) + \frac{9}{2} \tau_{2,24}(n),\\ 
{\mathcal N}_2(1,16; n) & = \frac{3}{2} \sigma(n) - \frac{9}{2} \sigma(n/2) - \frac{9}{2} \sigma(n/3) + 9 \sigma(n/4) + \frac{27}{2} \sigma(n/6) - 18 \sigma(n/8)\\
& ~\quad  - 27 \sigma(n/12) + 24\sigma(n/16) + 54 \sigma(n/24) + 72 \sigma(n/48) + \frac{9}{2} \tau_{2,48}(n),
\end{split}
\end{equation}
where $\tau_{2,24}(n)$ and $\tau_{2,48}(n)$ are Fourier coefficients of certain eta-quotients defined in \eqref{newforms}.
\end{thm}

\begin{thm}\label{Q3}
Let $n\in \mathbb N$. For each entry $(a_1, a_2, b_1)$ corresponding to ${\mathcal Q}_3$ in Table 1, the associated theta series is a modular form of 
weight $2$ on $\Gamma_0(48)$ with character $\chi$. Therefore, using the basis given in Table B (in \S 3.5), we have 
\begin{equation}
{\mathcal N}_3(a_1, a_2, b_1; n) = \sum_{i=1}^{\ell_\chi} \beta_{i,\chi} A_{i,\chi}(n),
\end{equation}
where $A_{i,\chi}(n)$ are the Fourier coefficients of the basis elements $f_{i,\chi}$ and the values of the constants $\beta_{i,\chi}$'s are given in (\S 6, Table 3).
\end{thm}

\section{Proofs} 
Here again, we shall take $\chi$ to be one of the four characters $\chi_0, \chi_8, \chi_{12}$ or $\chi_{24}$ and $\ell_\chi$ is the dimension of the 
space of modular forms $M_2(48, \chi)$. The main ingredient in proving our theorems is the construction of explicit bases for the spaces $M_2(48,\chi)$.
For uniformity, we shall denote these basis elements as $\{f_{i,\chi}(z): 1\le i\le \ell_\chi\}$ and write their Fourier expansions as 
\begin{equation}
f_{i,\chi}(z) = \sum_{n\ge 0} A_{i,\chi}(n) e^{2\pi inz}.
\end{equation}
The basis elements $f_{i,\chi}(z)$ are explicitly given in \S 3.5.

\subsection{A basis for $M_2(48, \chi_0)$}
The vector space $M_2(48, \chi_0)$ has dimension $14$ with \\
$\dim_{\mathbb C} {\mathcal E}_2(48, \chi_0)=11$ and $\dim_\mathbb C{S_2(48, \chi_0)}=3$. 
For $a, b$ divisors of $N$ with $a|b, (b>a)$ we define $\phi_{a, b}(z)$ to be 
\begin{equation}\label{phi-ab}
\phi_{a, b}(z) := \frac{1}{b-a}(bE_2(bz)-aE_2(az)).
\end{equation} 
It is easy to see that $\phi_{a, b}\in M_2(N, \chi_0).$ We need the following two eta-quotients:  
\begin{equation}\label{newforms}
\Delta_{2, 24}(z) = 2^14^16^112^1=\sum_{n=1}^\infty \tau_{2, 24}(n)q^n, \quad 
\Delta_{2, 48}(z) = 2^{-1}4^4 6^{-1}8^{-1}12^424^{-1}
=\sum_{n=1}^\infty \tau_{2, 48}(n)q^n.
\end{equation}
Using the above functions, we give a basis for the space $M_2(48, \chi_0)$  in the following proposition.

\begin{prop}\label{chi0}
A basis for the space of Eisenstein series ${\mathcal E}_2(48, \chi_0)$ is given by 
$$
\{ \phi_{1, b}: b|48 (b>1), E_{2, \chi_{-4}, \chi_{-4}}(z), E_{2, \chi_{-4}, \chi_{-4}}(3z)\} 
$$
and a basis for the space of cusp forms ${S}_2(48, \chi_0)$ is given by 
$$\{\Delta_{2, 24}(z), \Delta_{2, 24}(2z), \Delta_{2, 48}(z)\}.$$
\end{prop}

\subsection{A basis for $M_2(48, \chi_8)$} 
The vector space $M_2(48, \chi_8)$ has dimension $12$ with 
\linebreak
$\dim_{\mathbb C} {\mathcal E}_2(48, \chi_0)=8$ and $\dim_\mathbb C{S_2(48, \chi_0)}=4$. For the space of cusp forms, we need the following eta-quotients.  
\begin{equation*}
\begin{split}
\Delta_{2, 24, \chi_8;1} (z) & =1^12^{-1}3^{-1}6^48^212^{-1} =\sum_{n=1}^\infty \tau_{2, 24, \chi_8;1}(n)q^n,\\ 
\Delta_{2, 24, \chi_8;2} (z) & =1^24^{-1}6^{-1}8^112^424^{-1}= \sum_{n=1}^\infty \tau_{2, 24, \chi_8;2}(n)q^n.
\end{split}
\end{equation*}
The following proposition gives a basis of the space $M_2(48, \chi_8)$.

\begin{prop}\label{chi8}
A basis for the space of Eisenstein series ${\mathcal E}_2(48, \chi_8)$ is given by 
$$
\{E_{2, {\bf 1}, \chi_{8}}(az), a|6,  E_{2, \chi_{8}, {\bf 1}}(bz), b|6\}
$$
and a basis for the space of cusp forms ${S}_2(48, \chi_8)$ is given by 
$$\{\Delta_{2, 24, \chi_8;1}(z), \Delta_{2, 24, \chi_8;1}(2z),  \Delta_{2, 24, \chi_8;2}(z), \Delta_{2, 24, \chi_8;2}(2z)\}.$$
\end{prop}

\subsection{A basis for $M_2(48, \chi_{12})$}
The vector space $M_2(48, \chi_{12})$ has dimension $14$  with 
\linebreak 
$\dim_\mathbb C{S_2(48, \chi_{12})}=2.$
For the space of cusp forms, we use the following eta-quotient. 
$$
\Delta_{2, 48, \chi_{12}}(z) =1^{-4}2^{11}4^{-5}6^18^112^{-1}24^1 = \sum_{n=1}^\infty a_{2, 48,\chi_{12}}(n) q^n.
$$
Using the above eta-quotient we define the following two cusp forms (which are obtained by considering the character twists of the above eta-quotient). 
\begin{equation}
\Delta_{2, 48, \chi_{12};1}(z) = \sum_{n\ge 1 \atop{n\equiv 1 \!\!\!\pmod{4}}} a_{2, 48,\chi_{12}}(n) q^n,
\quad \Delta_{2, 48, \chi_{12};2}(z) = \sum_{n\ge 1 \atop{n\equiv 3 \!\!\!\pmod{4}}} a_{2, 48,\chi_{12}}(n) q^n.
\end{equation}
Using these functions we give a basis for the space $M_2(48, \chi_{12})$ in the following proposition.

\smallskip

\begin{prop}\label{basischi12}
A basis for the space of Eisenstein series ${\mathcal E}_2(48, \chi_{12})$ is given by 
$$
\{E_{2,\bf{1}, \chi_{12}}(az); a|4,  E_{2, \chi_{12}, \bf{1}}(bz); b|4, E_{2, \chi_{-4}, \chi_{-3}}(t_1z); t_1|4, E_{2, \chi_{-3}, \chi_{-4}}(t_2z); t_2|4\}
$$
and a basis for the space of cusp forms ${S}_2(48, \chi_{12})$ is given by 
$$\{\Delta_{2, 48, \chi_{12};1}(z), \Delta_{2, 48, \chi_{12};2}(z)\}.$$
\end{prop}

\smallskip

\subsection{A basis for $M_2(48, \chi_{24})$} 
The vector space $M_2(48, \chi_{24})$ has dimension $12$ with 
\linebreak
$\dim_\mathbb C{S_2(48, \chi_{12})}=4$. We need the following eta-quotients.
\begin{equation*}
\begin{split}
\Delta_{2, 24, \chi_{24};1}(z) & = 1^12^{-1}3^{-1}4^16^412^{-2}24^2 = \sum_{n=1}^\infty \tau_{2, 24, \chi_{24}; 1} q^n, \\
\Delta_{2, 48, \chi_{24};2}(z) & =1^22^{-2}4^46^18^{-1}12^{-1}24^1 = \sum_{n=1}^\infty \tau_{2, 24, \chi_{24}; 2} q^n.
\end{split}
\end{equation*}
In the following we give a basis for the space $M_2(48, \chi_{24})$.

\smallskip

\begin{prop}\label{basischi24}
A basis for the space of Eisenstein series ${\mathcal E}_2(48, \chi_{24})$ is given by 
$$
\{E_{2, {\bf 1}, \chi_{24}}(az); a|2,  E_{2, \chi_{24}, {\bf 1}}(bz); b|2, E_{2, \chi_{-3}, \chi_{-8}}(t_1z); t_1|2, E_{2, \chi_{-8}, \chi_{-3}}(t_2z); t_2|2\}
$$
and a basis for the space of cusp forms ${S}_2(48, \chi_{24})$ is given by 
$$
\{\Delta_{2, 24, \chi_{24};1}(z),\Delta_{2, 24, \chi_{24};1}(2z), \Delta_{2, 24, \chi_{24};2}(z), \Delta_{2, 24, \chi_{24};2}(2z)\}.
$$
\end{prop}

\smallskip

\subsection{Combined table for bases}
In this section we combine all the bases given in Propositions 3.1 to 3.4 in a tabular form along with identifying the elements $f_{i,\chi}(z)$ for each $i$, $1\le i\le \ell_\chi$. 

\bigskip

\smallskip

\begin{center}
Table B (List of basis elements)\\
\begin{tabular}{|lll|}
\hline
$f_{1,\chi_0}(z) =  \phi_{1, 2}(z)$,& $f_{6,\chi_0}(z) =  \phi_{1, 12}(z)$, & $f_{11,\chi_0}(z) = E_{2, \chi_{-4}, \chi_{-4}}(3z)$, \\
$f_{2,\chi_0}(z) =  \phi_{1, 3}(z)$, & $f_{7,\chi_0}(z) =  \phi_{1, 16}(z)$,&  $f_{12,\chi_0}(z) = \Delta_{2, 24}(z)$,\\
 $f_{3,\chi_0}(z)=  \phi_{1, 4}(z)$, &  $f_{8,\chi_0}(z) =  \phi_{1, 24}(z)$, &  $f_{13,\chi_0}(z) = \Delta_{2, 24}(2z)$, \\
  $f_{4,\chi_0}(z) =  \phi_{1, 6}(z)$, & $f_{9,\chi_0}(z) =  \phi_{1, 48}(z)$, &  $f_{14,\chi_0}(z) = \Delta_{2, 48}(z).$\\
 $f_{5,\chi_0}(z) =  \phi_{1, 8}(z)$, &   $f_{10,\chi_0}(z) = E_{2, \chi_{-4}, \chi_{-4}}(z)$,&\\
\hline
\hline
$f_{1,\chi_8}(z) = E_{2, {\bf 1}, \chi_{8}}(z)$, & $f_{5,\chi_8}(z)=E_{2, \chi_{8}, {\bf 1}}(z)$, & $f_{9,\chi_8}(z)=\Delta_{2, 24, \chi_8;1}(z),$\\
$f_{2,\chi_8}(z) = E_{2, {\bf 1}, \chi_{8}}(2z)$, & $ f_{6,\chi_8}(z)=E_{2, \chi_{8}, {\bf 1}}(2z)$, & $f_{10,\chi_8}(z)=\Delta_{2, 24, \chi_8;1}(2z),$\\
$f_{3,\chi_8}(z) = E_{2, {\bf 1}, \chi_{8}}(3z)$, & $f_{7,\chi_8}(z)=E_{2, \chi_{8}, {\bf 1}}(3z)$, & $  f_{11,\chi_8}(z)=\Delta_{2, 24, \chi_8;2}(z),$\\
$f_{4,\chi_8}(z) = E_{2, {\bf 1}, \chi_{8}}(6z)$, & $ f_{8,\chi_8}(z)=E_{2, \chi_{8}, {\bf 1}}(6z)$, & $ f_{12,\chi_8}(z)=\Delta_{2, 24, \chi_8;2}(2z).$ \\
\hline
\hline
$f_{1,\chi_{12}}(z) =  E_{2, \textbf{1}, \chi_{12}}(z),$&$ f_{6,\chi_{12}}(z) = E_{2, \chi_{12},\textbf{1}}(4z),$ & $ f_{11,\chi_{12}}(z) = E_{2, \chi_{-3}, \chi_{-4}}(2z),$\\
$f_{2,\chi_{12}}(z)) = E_{2, \textbf{1}, \chi_{12}}(2z),$& $f_{7,\chi_{12}}(z) = E_{2, \chi_{-4}, \chi_{-3}}(z),$& $f_{12,\chi_{12}}(z) = E_{2, \chi_{-3}, \chi_{-4}}(4z),$\\
$f_{3,\chi_{12}}(z) = E_{2, \textbf{1}, \chi_{12}}(4z), $ &$f_{8,\chi_{12}}(z) = E_{2, \chi_{-4}, \chi_{-3}}(2z),$& $f_{13,\chi_{12}}(z) = \Delta_{2, 48, \chi_{12};1}(z),$\\
$f_{4,\chi_{12}}(z) = E_{2, \chi_{12}, \textbf{1}}(z),$ & $ f_{9,\chi_{12}}(z) = E_{2, \chi_{-4}, \chi_{-3}}(4z),$ & $ f_{14,\chi_{12}}(z) = \Delta_{2, 48, \chi_{12};2}(z).$\\
$f_{5,\chi_{12}}(z)  = E_{2, \chi_{12}, \textbf{1}}(2z),$ &$f_{10,\chi_{12}}(z) = E_{2, \chi_{-3}, \chi_{-4}}(z),$ &\\
\hline
\hline
$f_{1,\chi_{24}}(z)=E_{2, \textbf{1}, \chi_{24}}(z),$ & $f_{5,\chi_{24}}(z)=E_{2, \chi_{-3}, \chi_{-8}}(z),$ & $f_{9,\chi_{24}}(z)=\Delta_{2,24, \chi_{24};1}(z),$\\
$ f_{2,\chi_{24}}(z)=E_{2, \textbf{1}, \chi_{24}}(2z),$& $f_{6,\chi_{24}}(z)=E_{2, \chi_{-3}, \chi_{-8}}(2z),$ & $ f_{10,\chi_{24}}(z)=\Delta_{2,24, \chi_{24};1}(2z),$\\
$ f_{3,\chi_{24}}(z)=E_{2, \chi_{24},\textbf{1}}(z), $ &$ f_{7,\chi_{24}}(z)=E_{2, \chi_{-8}, \chi_{-3}}(z), $ & $f_{11,\chi_{24}}(z)= \Delta_{2,24, \chi_{24};2}(z),$\\
$ f_{4,\chi_{24}}(z)=E_{2, \chi_{24},\textbf{1}}(2z),$ & $ f_{8,\chi_{24}}(z)=E_{2, \chi_{-8}, \chi_{-3}}(2z), $& $ f_{12,\chi_{24}}(z)=\Delta_{2,24, \chi_{24};2}(2z).$\\
\hline
\end{tabular}
\end{center}


\bigskip

We are now ready to prove the theorems. The generating functions for the two types of quadratic forms considered in this paper, viz., sum of squares and forms of type 
$x^2+xy+y^2$ are given respectively by the classical theta function 
\begin{equation}\label{theta}
\Theta(z) = \sum_{n\in {\mathbb Z}} e^{2\pi i n^2z},
\end{equation}
and the function 
\begin{equation}\label{F}
{\mathcal F}(z) = \sum_{m,n\in {\mathbb Z}} e^{2\pi i (m^2+mn+n^2)z}.
\end{equation}
The theta function $\Theta(z)$ is a modular form of weight $1/2$ on $\Gamma_0(4)$ and ${\mathcal F}(z)$ is a modular form of weight $1$ on $\Gamma_0(3)$ with character $\left(\frac{\cdot}{3}\right)$ (see \cite{koblitz}, \cite[Theorem 4]{schoeneberg}, \cite{borwein} for details). To each quadratic form $(a_1, a_2, a_3, a_4)$ as in the Table 1 (corresponding to 
the quadratic forms ${\mathcal Q}_1$), the associated theta series is given by 
\begin{equation}\label{theta-q1}
\Theta(a_1z) \Theta(a_2z) \Theta(a_3z) \Theta(a_4z).
\end{equation}
By using \cite[Lemmas 1--3]{alladi-proc}, we see that the above function is a modular form in $M_2(48, \chi)$, where $\chi$ is one of the four characters that appear in Table 1. Now using the bases constructed as in Table B, one can express each of the theta products \eqref{theta-q1} as a linear combination of the respective basis elements. Since ${\mathcal N}_1(a_1,a_2,a_3,a_4;n)$ is the $n$-th Fourier coefficient of the theta product \eqref{theta-q1}, by comparing the $n$-th Fourier coefficients, we get the required formulae in \thmref{Q1}. 

\smallskip

Next, for the four quadratic forms given by the pairs $(1,2), (1,4), (1,8), (1,16)$ in Table 1, the corresponding theta series is the product of the forms ${\mathcal F}(b_1z)$ and ${\mathcal F}(b_2z)$. Again by using Lemmas 1 and 3 in \cite{alladi-proc}, these forms belong $M_2(48, \chi_0)$. So, we can express these 4 forms as a linear 
combination of the basis elements of $M_2(48, \chi_0)$, which we denote as follows. Let $(b_1,b_2) \in \{(1,2), (1,4), (1,8), (1,16)\}$. Then 
\begin{equation}
{\mathcal N}_2(b_1,b_2;n) = \sum_{i=1}^{14} c_i A_{i,\chi_0}(n),
\end{equation}
where $A_{i,\chi_0}(n)$ are the Fourier coefficients of the basis elements $f_{i,\chi_0}(z)$ (given in Table B). The values of the constants $c_i$ for each pair $(b_1,b_2)$ are given in the following table. \\



\begin{center}
{\tiny
{Table C}\\
\smallskip
\begin{tabular}{|c|c|c|c|c|c|c|c|c|c|c|c|c|c|c|}
\hline
$b_1,b_2$&$c_1$ &$c_2$& $c_3$& $c_4$&$c_5$ &$c_6$& $c_7$& $c_8$&$c_9$ &$c_{10}$&$c_{11}$ &$c_{12}$& $c_{13}$& $c_{14}$\\
\hline 
&&&&&&&&&&&&&&\\ 
$1,2$&$\frac{1}{4}$&$\frac{-1}{2}$&0&$\frac{5}{4}$&0&0&0&0&0&0&0&0&0&0\\
&&&&&&&&&&&&&&\\ 
$1,4$&$\frac{3}{8}$&$\frac{1}{2}$&$\frac{-3}{4}$&$\frac{-15}{8}$&0&$\frac{11}{4}$&0&0&0&0&0&0&0&0\\
&&&&&&&&&&&&&&\\ 
$1,8$&$\frac{3}{32}$&$\frac{-1}{8}$&$\frac{-9}{32}$&$\frac{15}{32}$&$\frac{7}{16}$&$\frac{-33}{32}$&0&$\frac{23}{16}$&0&0&0&$\frac{9}{2}$&0&0\\
&&&&&&&&&&&&&&\\ 
$1,16$&$\frac{3}{32}$&$\frac{1}{8}$&$\frac{-9}{32}$&$\frac{-15}{32}$&$\frac{21}{32}$&$\frac{33}{32}$&$\frac{-15}{16}$&$\frac{-69}{32}$&$\frac{47}{16}$&0&0&0&0&$\frac{9}{2}$\\
&&&&&&&&&&&&&&\\ 
\hline
\end{tabular}
}
\end{center} 


The values of $c_i$ are non-zero only in the case of basis elements which are either $\phi_{1,b}(z)$, $b\vert 48$ and $b>1$ or one of the cusp forms $\Delta_{2,24}(z)$, 
$\Delta_{2,48}(z)$. The Fourier expansion of the Eisenstein series $\phi_{a,b}(z)$ is given as follows. 
\begin{equation*}
\phi_{a,b}(z) = 1 + \frac{24 a}{b-a}\sum_{n\ge 1} \sigma(n/a) q^n - \frac{24 b}{b-a}\sum_{n\ge 1} \sigma(n/b) q^n.
\end{equation*}
By substituting the values of the constants $c_i$ in the expression along with the Fourier expansion of the above basis elements, we get the required formulas in 
\thmref{Q2}. 

\smallskip

Finally, the theta series corresponding to each quadratic form ${\mathcal Q}_3$ represented by the triplets $(a_1,a_2,b_1)$ in Table 1 is the product 
$\Theta(a_1 z) \Theta(a_2 z) {\mathcal F}(b_1 z)$. By using Lemmas 1 to 3 of \cite{alladi-proc}, this is a modular form of weight $2$ on $\Gamma_0(48)$ with 
one of the characters $\chi_0$ or $\chi_d$, $d=8,12,24$ (depending on the triplets $(a_1,a_2,b_1)$). Formulas in \thmref{Q3} now follow from comparing the Fourier 
coefficients of these associated modular forms. 

This completes the proofs of the theorems.

\section{Sample Formulae}
In this section we give explicit formulas for a few cases (5 examples for ${\mathcal Q}_1$ and 6 examples for ${\mathcal Q}_3$). 

\begin{eqnarray*}
{\mathcal N}_1(1, 2, 4, 4; n) &=& - 2 \sigma_{2, {\bf 1}, \chi_{8}}(n/2) +  2 \sigma_{2, \chi_{8}, {\bf 1}}(n),\\ 
{\mathcal N}_1(1, 2, 4, 6; n) &=& - \sigma_{2, {\bf 1}, \chi_{12}}(n/4) +  \frac{3}{2} \sigma_{2, \chi_{12}, {\bf 1}}(n) 
+  \frac{1}{2} \sigma_{2, \chi_{-4}, \chi_{-3}}\left({n} \right) -3  \sigma_{2, \chi_{-3}, \chi_{-4}}(n/4), \\
{\mathcal N}_1(1, 2, 4, 12)(n) &=& \sigma_{2,\chi_{24}, {\bf 1}} (n)-\frac{1}{3} \sigma_{2, {\bf 1}, \chi_{24}}(n/2) + \frac{1}{3} \sigma_{2,\chi_{-8}, \chi_{-3}}\left({n}\right) + \sigma_{2,\chi_{-3}, \chi_{-8}}(n/2)\\
&&+4 \tau_{2, 24, \chi_{24}; 1} (n/2) + \frac{2}{3} \tau_{2, 24, \chi_{24}; 2} \left({n} \right)+ \frac{4}{3} \tau_{2, 24, \chi_{24}; 2} (n/2),\\
&&\\
{\mathcal N}_1(1, 3, 4, 6;n) & = & -\frac{4}{5} \sigma_{2, {\bf 1}, \chi_8}(n/2) -\frac{6}{5}  \sigma_{2, {\bf 1}, \chi_8}(n/6) +\frac{8}{5} \sigma_{2, \chi_8, {\bf 1}}\left({n}\right) -\frac{12}{5} \sigma_{2, \chi_8, {\bf 1}}(n/3) \\
&&+ \frac{8}{5} \tau_{2, 24, \chi_8; 1}(n) - \frac{8}{5} \tau_{2, 24, \chi_8; 1} (n/2)-
\frac{6}{5} \tau_{2, 24, \chi_8; 2}(n) - \frac{8}{5} \tau_{2, 24, \chi_8; 2} (n/2),\\
&&\\
{\mathcal N}_1(1, 3, 4, 12; n) &=&\frac{1204}{1081} \sigma(n)-3\sigma(n/2)-3\sigma(n/3)+ 10 \sigma(n/4)+ 9 \sigma(n/6) -12 \sigma(n/8) \\
& & -30 \sigma(n/12)+ \frac{360}{23} \sigma(n/24)+ \frac{1656}{47} \sigma(n/48) 
- \frac{47}{24} \sigma_{2, \chi_{-4}, \chi_{-4}}(n)+ \tau_{2, 48}(n),\\
\mathcal{N}_3(1, 3, 1; n) &=&  8\sigma(n)-12\sigma(n/2)-24\sigma(n/3) +16\sigma(n/4) + 36\sigma(n/6)- 48 \sigma(n/12),\\
{\mathcal N}_3(1, 3, 16; n) &=&-\frac{17}{92} \sigma(n)-\frac{3}{2}\sigma(n/2)-\frac{3}{2}\sigma(n/3)+7\sigma(n/4)+ \frac{9}{2}\sigma(n/6)-18\sigma(n/8) \\
& &-21 \sigma(n/12)+24 \sigma\left(\frac{n}{16} \right)+ 54 \sigma(n/24)-72\sigma(n/48) + \frac{3}{2}\tau_{2, 48}(n),\\
{\mathcal N}_3(1, 4, 8; n) &=&  \frac{1}{4}\sigma_{2, {\bf 1}, \chi_{12}}(n)- \frac{1}{4}\sigma_{2,{\bf 1}, \chi_{12}}(n/2)-\sigma_{2, {\bf 1}, \chi_{12}}(n/4)+\frac{3}{4}
\sigma_{2, \chi_{12},{\bf 1}}(n)\\ 
&&-\frac{3}{2}\sigma_{2, \chi_{12}, {\bf 1}}(n/2)+6\sigma_{2, \chi_{12}, {\bf 1}}(n/4)+ \frac{1}{4} \sigma_{2,\chi_{-4}, \chi_{-3}}(n)+ \frac{1}{2} 
\sigma_{2,\chi_{-4}, \chi_{-3}}(n/2)\\
&& + 2 \sigma_{2,\chi_{-4}, \chi_{-3}}(n/4)+ \frac{3}{4} \sigma_{2,\chi_{-3}, \chi_{-4}}(n)+ \frac{3}{4} \sigma_{2,\chi_{-3}, \chi_{-4}}(n/2) - 3 \sigma_{2,\chi_{-3}, \chi_{-4}}(n/4),\\
{\mathcal N}_3(2, 3, 1;n) & = & \frac{2}{5} \sigma_{2,{\bf 1}, \chi_8}(n)-\frac{12}{5} \sigma_{2,{\bf 1}, \chi_8}(n/3)+\frac{16}{5} \sigma_{2, \chi_8,{\bf 1}}(n) + 
\frac{96}{5} \sigma_{2, \chi_8, {\bf 1}}(n/3)\\
& & +\frac{12}{5} \tau_{2, 24,\chi_8; 2}(n/3),\\
\end{eqnarray*}
\begin{eqnarray*}
\mathcal{N}_3(3, 3, 4; n) &=&  - \sigma_{2,{\bf 1}, \chi_{12}}(n)+\sigma_{2, \chi_{12},{\bf 1}}(n)- \sigma_{2,\chi_{-4}, \chi_{-3}}(n)+\sigma_{2,\chi_{-3}, \chi_{-4}}(n),\\
{\mathcal N}_3(3, 6, 2;n) & = & -\frac{1}{3} \sigma_{2,1, \chi_{24}}(n)+\frac{4}{3} \sigma_{2, \chi_{24},1}(n)-\frac{4}{3} \sigma_{2, \chi_{-3}, \chi_{-8}}(n)+\frac{1}{3} 
\sigma_{2, \chi_{-8}, \chi_{-3}}(n)\\
& & +\frac{4}{3} \tau_{2, 24,\chi_{24}; 1}(n).\\
\end{eqnarray*}


\section{Remarks on equivalence of formulas} 
As mentioned in the introduction, there are 36 formulas among our results which are known already using different methods (19 corresponding to ${\mathcal Q}_1$, 
2 corresponding to ${\mathcal Q}_2$, 15 corresponding to ${\mathcal Q}_3$). The two cases $(1,2), (1,4)$ corresponding to ${\mathcal Q}_2$ are the same formulas. 
However, for the remaining 34 cases, the earlier formulas had different expressions. In this section, we shall show that our formulas are equivalent to these 
earlier results. Here we demonstrate three cases ${\mathcal N}_1(1,2,4,4;n)$,  ${\mathcal N}_3(3,3,1;n)$ and ${\mathcal N}_3(3,3,4;n)$, which appear in \S 4. 
The remaining cases follow using similar arguments. 

We start with our expressions and derive the previous formulas. First, we consider the case $(1,2,4,4)$.
\begin{equation*}
{\mathcal N}_1(1, 2, 4, 4; n)  ~=~ - 2 \sigma_{2, {\bf 1}, \chi_{8}}(n/2) +  2 \sigma_{2, \chi_{8}, {\bf 1}}(n) ~=~ - 2 R(n/2) + 2 S(n),\\
\end{equation*}
\begin{equation*}
{\rm where} \quad \qquad S(n)  = \sum_{d\vert n} \left(\frac{8}{n/d}\right) d, \quad 
R(n)  = \sum_{d\vert n} \left(\frac{8}{d}\right) d. \hskip 2cm \mbox{(\cite[Eqs.(4.1),(4.2)]{aalw3})}
\end{equation*}
When $n= 2^\alpha N$, $\alpha\ge 0$ and $N$ is odd, we have $R(n) = \left(\frac{8}{N}\right) S(N)$ and $S(n) = 2^\alpha S(N)$. Now, following the simplifications as 
done in \cite[p. 31--32]{aalw3}, we see that 
\begin{equation}
{\mathcal N}_1(1, 2, 4, 4; n) = \left(2^{\alpha +1} - (1+(-1)^n)\left(\frac{8}{N}\right)\right) S(N),
\end{equation}
which is Theorem 5.4 of \cite{aalw3}.  
Next, we consider the case  $(1,3,1)$. 
\begin{equation*}
\mathcal{N}_3(1, 3, 1; n) =  8\sigma(n)-12\sigma(n/2)-24\sigma(n/3) +16\sigma(n/4) + 36\sigma(n/6)- 48 \sigma(n/12).\\
\end{equation*}
Write $n\in {\mathbb N}$ as $n= 2^\alpha 3^\beta N$, where $\gcd(N,6)=1$. If $n\equiv 1\pmod{2}$, then all the terms except for the first and third terms vanish.
Therefore, 
\begin{equation*}
\mathcal{N}_3(1, 3, 1; n) ~=~ 8 \sigma(3^\beta N) - 24 \sigma(3^{\beta-1}N) ~=~ 8 \sigma(N) \left(\sigma(3^\beta) - 3 \sigma(3^{\beta-1})\right)
~ =~ 8 \sigma(N).
\end{equation*}
When $n$ is even, a similar calculation leads to 
\begin{equation*}
\begin{split}
\mathcal{N}_3(1, 3, 1; n) &= 8 \sigma(2^\alpha 3^\beta N) - 12 \sigma(2^{\alpha-1}3^\beta N) -24 \sigma(2^\alpha 3^{\beta-1}N) + 16  \sigma(2^{\alpha-2}3^\beta N) \\
& \quad + 36  \sigma(2^{\alpha-1}3^{\beta-1} N) - 48  \sigma(2^{\alpha-2}3^{\beta -1} N)\\   
& = 4\sigma(N) \left[\sigma(3^\beta)\left(2\sigma(2^{\alpha}) - 3 \sigma(2^{\alpha-1}) + 4 \sigma(2^{\alpha-2})\right)\right.\\
& \quad \left. + 3 \sigma(3^{\beta-1})
\left(-2\sigma(2^{\alpha}) + 3 \sigma(2^{\alpha-1}) - 4 \sigma(2^{\alpha-2})\right)\right]\\
& = 12\sigma(N) (2^\alpha -1)\left(\sigma(3^\beta) - 3\sigma(3^{\beta-1})\right) ~=~ 12 (2^\alpha -1) \sigma(N).\\
\end{split}
\end{equation*}
The above two cases are the formulas that appear in the first formula in Theorem 1.2 (iii) of \cite[p. 1400]{apw}. 
 Finally we consider the formula for ${\mathcal N}_3(3,3,4;n)$. Let $A(n), B(n), C(n), D(n)$ be as in Definition 3.1 of \cite{aalw2}. Then 
\begin{equation*}
\begin{split}
\mathcal{N}_3(3, 3, 4; n) & =  - \sigma_{2,{\bf 1}, \chi_{12}}(n)+\sigma_{2, \chi_{12},{\bf 1}}(n)- \sigma_{2,\chi_{-4}, \chi_{-3}}(n)+\sigma_{2,\chi_{-3}, \chi_{-4}}(n),\\
& = -D(n) + A(n) - C(n) + B(n),\\
\end{split}
\end{equation*}
from which the formula in Theorem 16.1 of \cite{aalw2} follows in the same way it was proved there.


\section{Tables for Theorems \ref{Q1} and \ref{Q3}}
In this section, we shall give Tables 2 and 3, which give explicit coefficients $\alpha_{i,\chi}$ and $\beta_{i,\chi}$ that appear in \thmref{Q1} and \thmref{Q3}.


\begin{center}
{\tiny
\textbf{Table 2 for the character} $\chi_{0}$.\\
\begin{tabular}{|c|c|c|c|c|c|c|c|c|c|c|c|c|c|c|} 
\hline
$\!\!\!a_1, a_2, a_3, a_4 \!\!\!$ &$\!\!\alpha_{1,\chi_{0}}\!\!$ & $\!\!\alpha_{2,\chi_{0}}\!\!$ & $\!\!\alpha_{3,\chi_{0}}\!\!$ & $\!\!\alpha_{4,\chi_{0}}\!\!$ & $\!\!\alpha_{5,\chi_{0}}\!\!$ & $\!\!\alpha_{6,\chi_{0}}\!\!$ & $ \!\!\alpha_{7,\chi_{0}}\!\!$ & $\!\!\alpha_{8,\chi_{0}}\!\!$ & $\!\!\alpha_{9,\chi_{0}}\!\!$ & $\!\!\alpha_{10,\chi_{0}}\!\!$ & $\!\!\alpha_{11,\chi_{0}}\!\!$ & $\!\!\alpha_{12,\chi_{0}}\!\!$ & $\!\!\alpha_{13,\chi_{0}}\!\!$ & $\!\!\alpha_{14,\chi_{0}}\!\!$\\ 
\hline
&&&&&&&&&&&&&&\\
1, 1, 1, 4 & 0 & 0& $\frac{5}{8}$ & 0& $-\frac{7}{8}$ & 0 & $\frac{5}{4}$ & 0 & 0& 2& 0& 0& 0& 0 \\ 
&&&&&&&&&&&&&&\\ 
1, 1, 4, 4 & $\frac{1}{24}$ & 0 & 0 & 0& $-\frac{7}{24}$ & 0 & $\frac{5}{4}$ & 0 & 0 & 2 & 0 & 0 & 0& 0\\
&&&&&&&&&&&&&&\\ 
1, 1, 3, 12 &$ \frac{1}{12}$ & $\frac{1}{6}$ & $-\frac{5}{16}$ & $-\frac{5}{12}$ & $\frac{7}{16}$ & $\frac{55}{48}$ & $-\frac{5}{8}$
& $-\frac{23}{16}$ & $\frac{47}{24}$ & 1 & 3 & 0 & 2 & 1 \\
&&&&&&&&&&&&&&\\ 
1, 1, 12, 12 & $\frac{1}{48}$ & $\frac{1}{12}$ & 0& $-\frac{5}{48}$ & $\frac{7}{48}$ & 0 & $-\frac{5}{8}$ & $- \frac{23}{48}$  & $\frac{47}{24}$ & 1& 3 & 1 & 2& 1\\
&&&&&&&&&&&&&&\\ 
1, 2, 2, 4 & $ \frac{1}{24}$ & 0 & $0$ & 0 & $-\frac{7}{24}$ & $0$ & $\frac{5}{4}$
& $0$ & $0$ & 0 & 0 & 0 & 0 & 0 \\
&&&&&&&&&&&&&&\\ 
1, 2, 6, 12 & $\frac{1}{96}$ & $-\frac{1}{24}$ & 0 & $\frac{5}{96}$ & $-\frac{7}{96}$ & 0& $\frac{5}{16}$ & $-\frac{23}{96}$ & $\frac{47}{48}$ & 0 & 0 &$\frac{1}{2}$ & 1 & 1 \\
&&&&&&&&&&&&&&\\ 

1, 3, 3, 4 & $\frac{1}{12}$ & $\frac{1}{6}$ & $-\frac{5}{16}$ & $-\frac{5}{12}$ & $\frac{7}{16}$ & $\frac{55}{48}$ & $-\frac{5}{8}$ & $-\frac{23}{16}$ & $\frac{47}{24}$ & -1 & -3 & 0 & -2 & 1\\
&&&&&&&&&&&&&&\\ 
1, 3, 4, 12 & $\frac{1}{16}$ & $\frac{1}{12}$ & $-\frac{5}{16}$ & $-\frac{5}{16}$ & $\frac{7}{16}$ & $\frac{55}{48}$ & $-\frac{5}{8}$
& $-\frac{23}{16}$ & $\frac{47}{24}$ & 0 & 0 & 0& 0& 1 \\
&&&&&&&&&&&&&&\\ 
1, 4, 4, 4 & $\frac{1}{16}$ & 0 & $-\frac{5}{16}$ &  0& 0& 0& $\frac{5}{4}$ & 0 & 0 & 1 & 0 & 0 & 0 & 0 \\
&&&&&&&&&&&&&&\\ 
1, 4, 6, 6 & $\frac{1}{48}$ & $\frac{1}{12}$ & 0 & $-\frac{5}{48}$ & $\frac{7}{48}$ & 0 & $-\frac{5}{8}$ & $-\frac{23}{48}$ &
 $\frac{47}{24}$ & 0 & 0 & 1 & -2 & 0 \\
&&&&&&&&&&&&&&\\  
1, 4, 12, 12 & $\frac{1}{32}$ & $\frac{1}{24}$ & $-\frac{5}{32}$ & $-\frac{5}{32}$ & $\frac{7}{24}$ & $\frac{55}{96}$ & $-\frac{5}{8}$ & $-\frac{23}{24}$ & $\frac{47}{24}$ & $\frac{1}{2}$ & $\frac{3}{2}$ & $\frac{1}{2}$ & 0 & $\frac{1}{2}$\\
&&&&&&&&&&&&&&\\  
2, 2, 3, 12 & $\frac{1}{48}$ & $\frac{1}{12}$ & 0 & $-\frac{5}{48}$ & $\frac{7}{48}$ & 0 & $-\frac{5}{8}$ & $-\frac{23}{48}$ &
 $\frac{47}{24}$ & 0 & 0 &-1 & 2 & 0\\ 
&&&&&&&&&&&&&&\\  
2, 3, 4, 6 & $\frac{1}{96}$ &  $-\frac{1}{24}$ & 0 &  $\frac{5}{96}$ &  $-\frac{7}{96}$ & 0 &  $\frac{5}{16}$ &  $-\frac{23}{96}$ &  $\frac{47}{48}$ & 0 & 0 &  $\frac{1}{2}$ & 1 & -1\\
&&&&&&&&&&&&&&\\  
3, 3, 4, 4 &  $\frac{1}{48}$  &  $\frac{1}{12}$ & 0&  $-\frac{5}{48}$ &  $\frac{7}{48}$ & 0 &  $-\frac{5}{8}$ &  $-\frac{23}{48}$ &  $\frac{47}{24}$ & -1 & -3 & -1 & -2& 1 \\
&&&&&&&&&&&&&&\\  
3, 4, 4, 12 &  $\frac{1}{32}$ &  $\frac{1}{24}$ &  $-\frac{5}{32}$ &  $-\frac{5}{32}$ &  $\frac{7}{24}$ &  $\frac{55}{96}$ &  $-\frac{5}{8}$ &  $-\frac{23}{24}$ &  $\frac{47}{24}$ &  $-\frac{1}{2}$ &  $-\frac{3}{2}$ &  $-\frac{1}{2}$ & 0 &  $\frac{1}{2}$\\
 &&&&&&&&&&&&&&\\ 
 \hline
\end{tabular}
}
\end{center}

\begin{center}
{\tiny
\textbf{Table 2 for the character} $\chi_{8}$.\\
\begin{tabular}{|c|c|c|c|c|c|c|c|c|c|c|c|c|} 
\hline
$a_1, a_2, a_3, a_4$ & $\alpha_{1,\chi_8}$ & $\alpha_{2,\chi_8}$ & $\alpha_{3,\chi_8}$ & $\alpha_{4,\chi_8}$ & $\alpha_{5,\chi_8}$ & $\alpha_{6,\chi_8}$ & $ \alpha_{7,\chi_8}$ &
$\alpha_{8,\chi_8}$ & $\alpha_{9,\chi_8}$ & $\alpha_{10,\chi_8}$ & $\alpha_{11,\chi_8}$ & $\alpha_{12,\chi_8}$ \\ 
\hline
1, 1, 2, 4 & 0 & -2 & 0 & 0 & 4 & 0 & 0 & 0 & 0 & 0 & 0 & 0 \\
  &&&&&&&&&&&&\\ 
1, 1, 6, 12 & 0 & $-\frac{4}{5}$ & 0 & $-\frac{6}{5}$ & $\frac{8}{5}$ & 0 & $-\frac{12}{5}$ & 0 & $\frac{8}{5}$ & $\frac{32}{5}$ & $\frac{4}{5}$ & $-\frac{8}{5}$ \\
 &&&&&&&&&&&&\\  
1, 2, 4, 4 & 0 & -2 & 0 & 0 & 2 & 0 & 0 & 0 & 0 & 0 & 0 & 0 \\
  &&&&&&&&&&&&\\ 
1, 2, 3, 12 & 0 & $\frac{2}{5}$ & 0 & $-\frac{12}{5}$ & $\frac{4}{5}$ & 0 & $\frac{24}{5}$ & 0 & $\frac{4}{5}$&  $\frac{24}{5}$ & $\frac{2}{5}$ & $-\frac{16}{5}$\\
  &&&&&&&&&&&&\\ 
1, 2, 12, 12 & 0 & $\frac{2}{5}$ & 0 & $-\frac{12}{5}$ & $\frac{2}{5}$ & 0 & $\frac{12}{5}$ &0 & $\frac{12}{5}$ & $\frac{24}{5}$ &
$-\frac{4}{5}$ & $-\frac{16}{5}$ \\
 &&&&&&&&&&&&\\  
1, 3, 4, 6 & 0 & $-\frac{4}{5}$ & 0 & $-\frac{6}{5}$ & $\frac{8}{5}$ & 0 & $-\frac{12}{5}$ & 0 & $\frac{8}{5}$ & $-\frac{8}{5}$ & $-\frac{6}{5}$ &  $-\frac{8}{5}$ \\
  &&&&&&&&&&&&\\ 
1, 4, 6, 12 & 0 & $-\frac{4}{5}$ & 0 & $-\frac{6}{5}$ & $\frac{4}{5}$ & 0 & $-\frac{6}{5}$ & 0 & $\frac{4}{5}$ & $\frac{12}{5}$ & $\frac{2}{5}$ & $-\frac{8}{5}$\\
 &&&&&&&&&&&&\\ 
2, 3, 3, 4 & 0 & $\frac{2}{5}$ & 0 & $-\frac{12}{5}$ & $\frac{4}{5}$ & 0 & $\frac{24}{5}$ & 0 & $-\frac{16}{5}$ & $-\frac{16}{5}$ & $\frac{12}{5}$ & $\frac{24}{5}$\\
 &&&&&&&&&&&&\\ 
2, 3, 4, 12 & 0 & $\frac{2}{5}$ & 0 & $-\frac{12}{5}$ & $\frac{2}{5}$ & 0 & $\frac{12}{5}$ & 0 & $-\frac{8}{5}$ & $\frac{4}{5}$ & $\frac{6}{5}$ & $\frac{4}{5}$\\
 &&&&&&&&&&&&\\ 
3, 4, 4, 6 & 0 & $-\frac{4}{5}$ & 0 & $-\frac{6}{5}$ & $\frac{4}{5}$ & 0 & $-\frac{6}{5}$ & 0 & $\frac{4}{5}$ & $-\frac{8}{5}$ & $-\frac{8}{5}$ & $-\frac{8}{5}$\\
\hline
\end{tabular}
}
\end{center}


\begin{center}
{\tiny
\textbf{Table 2 for the character} $\chi_{12}$.\\
\begin{tabular}{|c|c|c|c|c|c|c|c|c|c|c|c|c|c|c|} 
\hline
$\!\!\!\!a_1, a_2, a_3, a_4\!\!\!\!$ &$\!\!\alpha_{1,\chi_{12}}\!\!$ & $\!\!\alpha_{2,\chi_{12}}\!\!$ & $\!\!\alpha_{3,\chi_{12}}\!\!$ & $\!\!\alpha_{4,\chi_{12}}\!\!$ & $\!\!\alpha_{5,\chi_{12}}\!\!$ & $\!\!\alpha_{6,\chi_{12}}\!\!$ & $ \!\!\alpha_{7,\chi_{12}}\!\!$ & $\!\!\alpha_{8,\chi_{12}}\!\!$ & $\!\!\alpha_{9,\chi_{12}}\!\!$ & $\!\!\alpha_{10,\chi_{12}}\!\!$ & $\!\!\alpha_{11,\chi_{12}}\!\!$ & $\!\!\alpha_{12,\chi_{12}}\!\!$ & $\!\!\alpha_{13,\chi_{12}}\!\!$ & $\!\!\alpha_{14,\chi_{12}}\!\!$\\ 
\hline
1, 1, 1, 12 & $-\frac{1}{2}$ & $\frac{1}{2}$& -1 & 3 & 3 & -12& -1 & 1 & 4 & $\frac{3}{2}$ & $\frac{3}{2}$ & 3 & 3 & 1 \\
 &&&&&&&&&&&&&&\\ 
1, 1, 3, 4 & $-\frac{1}{2}$ & $\frac{1}{2}$& -1 & 3 & -3 & 12& -1 & -1 & -4 & $\frac{3}{2}$ & $\frac{3}{2}$ & 3 & 1 & -1 \\
 &&&&&&&&&&&&&&\\ 
1, 1, 4, 12 & $-\frac{1}{2}$ & $\frac{1}{2}$& -1 & $\frac{3}{2}$ & 0 & 0 & $-\frac{1}{2}$& 0 & 0 & $\frac{3}{2}$ & 
$\frac{3}{2}$& 3 & 2 & 0 \\
 &&&&&&&&&&&&&&\\ 
1, 2, 2, 12 & 0 & 0 & -1 & $\frac{3}{2}$ & 0 & 0 & $-\frac{1}{2}$& 0 & 0 & 0 & 0 & 3 & 1 & 1 \\
 &&&&&&&&&&&&&&\\ 
1, 2, 4, 6 & 0 & 0 & -1 & $\frac{3}{2}$ & 0 & 0 & $\frac{1}{2}$& 0 & 0 & 0 & 0 & -3 & 0 & 0 \\
 &&&&&&&&&&&&&&\\ 
1, 3, 3, 12 & $-\frac{1}{2}$ & $\frac{1}{2}$& -1 & 1 & -1 & 4 & 1 & 1 & 4 & $-\frac{1}{2}$ & $-\frac{1}{2}$& -1 & 1 & $\frac{1}{3}$\\
 &&&&&&&&&&&&&&\\ 
1, 3, 4, 4 & 0 & 0 & -1 & $\frac{3}{2}$ & -3& 12 &  $-\frac{1}{2}$& -1 & -4 & 0 & 0 & 3 & 1 & -1\\
 &&&&&&&&&&&&&&\\ 
1, 3, 12, 12 & 0 & 0 & -1 & $\frac{1}{2}$ & -1 & 4 &  $\frac{1}{2}$& 1 & 4 & 0 & 0 & -1 & 1 & $\frac{1}{3}$\\
 &&&&&&&&&&&&&&\\ 
1, 4, 4, 12 & $-\frac{1}{4}$ & $\frac{1}{4}$& -1 & $\frac{3}{4}$ & $-\frac{3}{2}$& 6 & $-\frac{1}{4}$ & $-\frac{1}{2}$& -2 & $\frac{3}{4}$ & $\frac{3}{4}$& 3 & 1 & 0\\
 &&&&&&&&&&&&&&\\ 
1, 6, 6, 12 & 0 & 0 & -1& $\frac{1}{2}$ &0 & 0 & $\frac{1}{2}$& 0 & 0 & 0& 0 &-1 & 1 & $-\frac{1}{3}$\\
 &&&&&&&&&&&&&&\\ 
1, 12, 12, 12 & $\frac{1}{4}$ & $-\frac{1}{4}$& 1 & $\frac{1}{4}$ & $-\frac{1}{2}$& 2& $\frac{1}{4}$ & $\frac{1}{2}$&2  & $\frac{1}{4}$& $\frac{1}{4}$ & -1 & 1 & 0 \\
 &&&&&&&&&&&&&&\\ 
2, 2, 3, 4 & 0 & 0 & -1 & $\frac{3}{2}$ & 0 & 0 & $-\frac{1}{2}$  & 0 & 0 & 0 & 0 & 3 & -1 & -1\\
 &&&&&&&&&&&&&&\\ 
2, 3, 6, 12 &  0 & 0 & -1 & $\frac{1}{2}$ & 0 & 0 & $-\frac{1}{2}$  & 0 & 0 & 0 & 0 & 1 & 0 & 0\\
 &&&&&&&&&&&&&&\\ 
3, 3, 3, 4 & $-\frac{1}{2}$ & $\frac{1}{2}$& -1 & 1 & 1 & -4& 1 & -1 & -4 & $-\frac{1}{2}$ & $-\frac{1}{2}$ & -1& -1& 1\\
 &&&&&&&&&&&&&&\\ 
3, 3, 4, 12 & $-\frac{1}{2}$ & $\frac{1}{2}$& -1 & $\frac{1}{2}$ & 0& 0 &$\frac{1}{2}$& 0 & 0 & $-\frac{1}{2}$ & $-\frac{1}{2}$& -1 & 0 & $\frac{2}{3}$\\
 &&&&&&&&&&&&&&\\ 
3, 4, 4, 4 & $\frac{1}{4}$ & $-\frac{1}{4}$ & -1 & $\frac{3}{4}$ & $-\frac{3}{2}$& 6 & $-\frac{1}{4}$ & $-\frac{1}{2}$& -2 & $-\frac{3}{4}$ & $-\frac{3}{4}$& 3 & 0 & -1\\
 &&&&&&&&&&&&&&\\ 
3, 4, 6, 6 & 0 & 0 & -1 & $\frac{1}{2}$ & 0 & 0 & $\frac{1}{2}$& 0 & 0 & 0& 0& -1 & -1 & $\frac{1}{3}$\\
 &&&&&&&&&&&&&&\\ 
3, 4, 12, 12 & $-\frac{1}{4}$ & $\frac{1}{4}$& -1 & $\frac{1}{4}$ & $-\frac{1}{2}$& 2 & $\frac{1}{4}$ & $\frac{1}{2}$& 2 & $-\frac{1}{4}$ & $-\frac{1}{4}$& -1 & 0 & $\frac{1}{3}$\\
\hline
\end{tabular}
}
\end{center}


\begin{center}
{\tiny
\textbf{Table 2 for the character} $\chi_{24}$.\\
\begin{tabular}{|c|c|c|c|c|c|c|c|c|c|c|c|c|c|cl} 
\hline
$\!\!a_1, a_2, a_3, a_4 \!\!$ & $\!\alpha_{1,\chi_{24}}\!$ & $\!\alpha_{2,\chi_{24}}\!$ & $\!\alpha_{3,\chi_{24}}\!$ & $\!\alpha_{4,\chi_{24}}\!$ & $\!\alpha_{5,\chi_{24}}\!$ & $\!\alpha_{6,\chi_{24}}\!$ & $\! \alpha_{7,\chi_{24}}\!$ & $\!\alpha_{8,\chi_{24}}\!$ & $\!\alpha_{9,\chi_{24}}\!$ & $\!\alpha_{10,\chi_{24}}\!$ & $\!\alpha_{11,\chi_{24}}\!$ & $\!\alpha_{12,\chi_{24}}\!$ \\ 
\hline
1, 1, 2, 12 & 0 & $-\frac{1}{3}$ &  2& 0 & $\frac{2}{3}$& 0 & 0 & 1 & 0 & 16 & $\frac{4}{3}$ & $\frac{16}{3}$\\
 &&&&&&&&&&&&\\ 
1, 1, 4, 6 & 0 & $-\frac{1}{3}$ &  2& 0 & $-\frac{2}{3}$& 0 & 0 & -1 & 8 & 0 & $\frac{8}{3}$ & $-\frac{8}{3}$\\
 &&&&&&&&&&&&\\ 
1, 2, 3, 4 & 0 & $-\frac{1}{3}$ &  2& 0 & $\frac{2}{3}$& 0 & 0 & 1 & 0 & -8 & $-\frac{2}{3}$ & $-\frac{8}{3}$\\
 &&&&&&&&&&&&\\ 
1, 2, 4, 12 & 0 & $-\frac{1}{3}$ &  1& 0 & $\frac{1}{3}$& 0 & 0 & 1 & 0 & 4 & $\frac{2}{3}$ & $\frac{4}{3}$\\
 &&&&&&&&&&&&\\ 
1, 3, 6, 12 & 0 & $-\frac{1}{3}$ & $\frac{2}{3}$& 0 & $\frac{2}{3}$ & 0 & 0 & $\frac{1}{3}$ & $\frac{4}{3}$ & $\frac{8}{3}$ & $\frac{2}{3}$ & 0 \\
 &&&&&&&&&&&&\\ 
1, 4, 4, 6 & 0 & $-\frac{1}{3}$ & 1 & 0 &  $-\frac{1}{3}$ & 0 & 0 & -1 & 4 & 0 & $\frac{4}{3}$ & $-\frac{8}{3}$ \\
 &&&&&&&&&&&&\\ 
1, 6, 12, 12 & 0 & $-\frac{1}{3}$ & $\frac{1}{3}$ & 0 & $\frac{1}{3}$& 0 & 0 & $\frac{1}{3}$ & $\frac{8}{3}$ & $\frac{8}{3}$ & $\frac{4}{3}$ & 0 \\
 &&&&&&&&&&&&\\ 
2, 3, 3, 12 &  0 & $-\frac{1}{3}$ & $\frac{2}{3}$ & 0 & $-\frac{2}{3}$& 0 & 0 & $-\frac{1}{3}$ & $-\frac{8}{3}$ & $\frac{16}{3}$ & 0&$\frac{8}{3}$ \\
 &&&&&&&&&&&&\\ 
2, 3, 4, 4 &   0 & $-\frac{1}{3}$ & 1 & 0 & $\frac{1}{3}$ &  0 & 0 &  1& 0 & -8 & $-\frac{4}{3}$ & $-\frac{8}{3}$ \\
 &&&&&&&&&&&&\\ 
2, 3, 12, 12 & 0 & $-\frac{1}{3}$ & $\frac{1}{3}$ &  0 &  $-\frac{1}{3}$ & 0& 0&  $-\frac{1}{3}$ &  $-\frac{4}{3}$ & $\frac{16}{3}$ & 0 & $\frac{8}{3}$\\
 &&&&&&&&&&&&\\ 
3, 3, 4, 6 & 0 &  $-\frac{1}{3}$ & $\frac{2}{3}$ &  0 &  $\frac{2}{3}$ & 0 & 0 & $\frac{1}{3}$ &  $-\frac{8}{3}$ & $-\frac{16}{3}$ & 
 $-\frac{4}{3}$ & 0 \\
  &&&&&&&&&&&&\\ 
3, 4, 6, 12 & 0 &  $-\frac{1}{3}$ & $\frac{1}{3}$ &  0 &  $\frac{1}{3}$ & 0 & 0 & $\frac{1}{3}$ &  $-\frac{4}{3}$ & $-\frac{4}{3}$ & 
 $-\frac{2}{3}$ & 0 \\
\hline
\end{tabular}
}
\end{center}

\smallskip

\newpage

\begin{center}
{\tiny
\textbf{Table 3 for the character} $\chi_{0}$.\\

\begin{tabular}{|c|c|c|c|c|c|c|c|c|c|c|c|c|c|c|}
\hline
&&&&&&&&&&&&&&\\
$a_1,a_2,b_1$&$\!\!\beta_{1,\chi_{0}}\!\!$ & $\!\!\beta_{2,\chi_{0}}\!\!$ & $\!\!\beta_{3,\chi_{0}}\!\!$ & $\!\!\beta_{4,\chi_{0}}\!\!$ & $\!\!\beta_{5,\chi_{0}}\!\!$ & $\!\!\beta_{6,\chi_{0}}\!\!$ & $ \!\!\beta_{7,\chi_{0}}\!\!$ & $\!\!\beta_{8,\chi_{0}}\!\!$ & $\!\!\beta_{9,\chi_{0}}\!\!$ & $\!\!\beta_{10,\chi_{0}}\!\!$ & $\!\!\beta_{11,\chi_{0}}\!\!$ & $\!\!\beta_{12,\chi_{0}}\!\!$ & $\!\!\beta_{13,\chi_{0}}\!\!$ & $\!\!\beta_{14,\chi_{0}}\!\!$\\ 
&&&&&&&&&&&&&&\\
\hline 
$1,3,1$&$\frac{1}{4}$&$\frac{2}{3}$&$\frac{-1}{2}$&$\frac{-5}{4}$&0&$\frac{11}{6}$&0&0&0&0&0&0&0&0\\
&&&&&&&&&&&&&&\\           
$1,3,2$&0&$\frac{-1}{6}$&$\frac{1}{4}$&0&0&$\frac{11}{12}$&0&0&0&0&0&0&0&0\\
&&&&&&&&&&&&&&\\           
$1,3,4$&$\frac{1}{8}$&$\frac{1}{6}$&$\frac{-1}{2}$&$\frac{-5}{8}$&0&$\frac{11}{6}$&0&0&0&0&0&0&0&0\\
&&&&&&&&&&&&&&\\
$1,3,8$&$\frac{1}{32}$&$\frac{-1}{24}$&$\frac{-7}{32}$&$\frac{5}{32}$&$\frac{7}{16}$&$\frac{-77}{96}$&0&$\frac{23}{16}$&0&0&0&$\frac{3}{2}$&0&0\\
&&&&&&&&&&&&&&\\
$1,3,16$&$\frac{1}{32}$&$\frac{1}{24}$&$\frac{-7}{32}$&$\frac{-5}{32}$&$\frac{21}{32}$&$\frac{77}{96}$&$\frac{-15}{16}$&$\frac{-69}{32}$&$\frac{47}{16}$&0&0&0&0&$\frac{3}{2}$\\
&&&&&&&&&&&&&&\\
$1,12,1$&$\frac{1}{8}$&$\frac{1}{3}$&$\frac{-5}{16}$&$\frac{-5}{8}$&$\frac{7}{16}$&$\frac{55}{48}$&$\frac{-5}{8}$&$\frac{-23}{16}$&$\frac{47}{24}$&1&3&0&6&3\\
&&&&&&&&&&&&&&\\
$1,12,2$&0&$\frac{-1}{12}$&$\frac{5}{32}$&0&$\frac{-7}{32}$&$\frac{55}{96}$&$\frac{5}{16}$&$\frac{-23}{32}$&$\frac{47}{48}$&$\frac{-1}{2}$&$\frac{3}{2}$&0&3&$\frac{3}{2}$\\
&&&&&&&&&&&&&&\\
$1,12,4$&$\frac{1}{16}$&$\frac{1}{12}$&$\frac{-5}{16}$&$\frac{-5}{16}$&$\frac{7}{16}$&$\frac{55}{48}$&$\frac{-5}{8}$&$\frac{-23}{16}$&$\frac{47}{24}$&1&3&0&0&0\\
&&&&&&&&&&&&&&\\
$1,12,8$&$\frac{1}{64}$&$\frac{-1}{48}$&$\frac{-5}{64}$&$\frac{5}{64}$&0&$\frac{-55}{192}$&$\frac{5}{16}$&0&$\frac{47}{48}$&$\frac{1}{4}$&$\frac{-3}{4}$&$\frac{3}{4}$&0&$\frac{3}{4}$\\
&&&&&&&&&&&&&&\\
$1,12,16$&$\frac{1}{64}$&$\frac{1}{48}$&$\frac{-5}{64}$&$\frac{-5}{64}$&$\frac{7}{32}$&$\frac{55}{192}$&$\frac{-5}{8}$&$\frac{-23}{32}$&$\frac{47}{24}$&$\frac{1}{4}$&$\frac{3}{4}$&$\frac{3}{4}$&0&$\frac{3}{4}$\\
&&&&&&&&&&&&&&\\
$2,6,1$&$\frac{7}{48}$&$\frac{-1}{4}$&$\frac{-3}{16}$&$\frac{35}{48}$&$\frac{7}{24}$&$\frac{-11}{16}$&0&$\frac{23}{24}$&0&0&0&3&0&0\\
&&&&&&&&&&&&&&\\
$3,4,1$&$\frac{1}{8}$&$\frac{1}{3}$&$\frac{-5}{16}$&$\frac{-5}{8}$&$\frac{7}{16}$&$\frac{55}{48}$&$\frac{-5}{8}$&$\frac{-23}{16}$&$\frac{47}{24}$&-1&-3&0&-6&3\\
&&&&&&&&&&&&&&\\
$3,4,2$&0&$\frac{-1}{12}$&$\frac{5}{32}$&0&$\frac{-7}{32}$&$\frac{55}{96}$&$\frac{5}{16}$&$\frac{-23}{32}$&$\frac{47}{48}$&$\frac{1}{2}$&$\frac{-3}{2}$&0&3&$\frac{-3}{2}$\\
&&&&&&&&&&&&&&\\
$3,4,4$&$\frac{1}{16}$&$\frac{1}{12}$&$\frac{-5}{16}$&$\frac{-5}{16}$&$\frac{7}{16}$&$\frac{55}{48}$&$\frac{-5}{8}$&$\frac{-23}{16}$&$\frac{47}{24}$&-1&-3&0&0&0\\
&&&&&&&&&&&&&&\\
$3,4,8$&$\frac{1}{64}$&$\frac{-1}{48}$&$\frac{-5}{64}$&$\frac{5}{64}$&0&$\frac{-55}{192}$&$\frac{5}{16}$&0&$\frac{47}{48}$&$\frac{-1}{4}$&$\frac{3}{4}$&$\frac{3}{4}$&0&$\frac{-3}{4}$\\
&&&&&&&&&&&&&&\\
$3,4,16$&$\frac{1}{64}$&$\frac{1}{48}$&$\frac{-5}{64}$&$\frac{-5}{64}$&$\frac{7}{32}$&$\frac{55}{192}$&$\frac{-5}{8}$&$\frac{-23}{32}$&$\frac{47}{24}$&$\frac{-1}{4}$&$\frac{-3}{4}$&$\frac{-3}{4}$&0&$\frac{3}{4}$\\
&&&&&&&&&&&&&&\\
$4,12,1$&$\frac{3}{16}$&$\frac{1}{4}$&$\frac{-7}{16}$&$\frac{-15}{16}$&$\frac{7}{16}$&$\frac{77}{48}$&$\frac{-5}{8}$&$\frac{-23}{16}$&$\frac{47}{24}$&0&0&0&0&3\\
\hline
\end{tabular}
}
\end{center} 

\begin{center}
{\tiny 
\textbf{Table 3 for the character} $\chi_{8}$.\\
\begin{tabular}{|c|c|c|c|c|c|c|c|c|c|c|c|c|}
\hline
&&&&&&&&&&&&\\
$a_1,a_2,b_1$& $\!\!\beta_{1,\chi_{8}}\!\!$ & $\!\!\beta_{2,\chi_{8}}\!\!$ & $\!\!\beta_{3,\chi_{8}}\!\!$ & $\!\!\beta_{4,\chi_{8}}\!\!$ & $\!\!\beta_{5,\chi_{8}}\!\!$ & $\!\!\beta_{6,\chi_{8}}\!\!$ & $ \!\!\beta_{7,\chi_{8}}\!\!$ & $\!\!\beta_{8,\chi_{8}}\!\!$ & $\!\!\beta_{9,\chi_{8}}\!\!$ & $\!\!\beta_{10,\chi_{8}}\!\!$ & $\!\!\beta_{11,\chi_{8}}\!\!$ & $\!\!\beta_{12,\chi_{8}}\!\!$\\ 
&&&&&&&&&&&&\\
\hline 
$1,6,1$&$\frac{-4}{5}$&0&$\frac{-6}{5}$&0&$\frac{32}{5}$&0&$\frac{-48}{5}$&0&$\frac{24}{5}$&0&$\frac{-12}{5}$&0\\
&&&&&&&&&&&&\\
$1,6,2$&$\frac{2}{5}$&0&$\frac{-12}{5}$&0&$\frac{8}{5}$&0&$\frac{48}{5}$&0&$\frac{12}{5}$&0&$\frac{-12}{5}$&0\\
&&&&&&&&&&&&\\
$1,6,4$&$\frac{-4}{5}$&0&$\frac{-6}{5}$&0&$\frac{8}{5}$&0&$\frac{-12}{5}$&0&0&0&$\frac{6}{5}$&0\\
&&&&&&&&&&&&\\
$1,6,8$&$\frac{2}{5}$&0&$\frac{-12}{5}$&0&$\frac{2}{5}$&0&$\frac{12}{5}$&0&$\frac{6}{5}$&0&0&0\\
&&&&&&&&&&&&\\
$1,6,16$&$\frac{2}{5}$&$\frac{-6}{5}$&$\frac{3}{5}$&$\frac{-9}{5}$&$\frac{2}{5}$&0&$\frac{-3}{5}$&0&$\frac{6}{5}$&$\frac{18}{5}$&0&$\frac{-12}{5}$\\
&&&&&&&&&&&&\\
$2,3,1$&$\frac{2}{5}$&0&$\frac{-12}{5}$&0&$\frac{16}{5}$&0&$\frac{96}{5}$&0&0&0&$\frac{12}{5}$&0\\
&&&&&&&&&&&&\\
$2,3,2$&$\frac{-4}{5}$&0&$\frac{-6}{5}$&0&$\frac{16}{5}$&0&$\frac{-24}{5}$&0&$\frac{-12}{5}$&0&0&0\\
&&&&&&&&&&&&\\
$2,3,4$&$\frac{2}{5}$&0&$\frac{-12}{5}$&0&$\frac{4}{5}$&0&$\frac{24}{5}$&0&$\frac{-12}{5}$&0&$\frac{6}{5}$&0\\
&&&&&&&&&&&&\\
$2,3,8$&$\frac{-4}{5}$&0&$\frac{-6}{5}$&0&$\frac{4}{5}$&0&$\frac{-6}{5}$&0&$\frac{6}{5}$&0&$\frac{-6}{5}$&0\\
&&&&&&&&&&&&\\
$2,3,16$&$\frac{-1}{5}$&$\frac{3}{5}$&$\frac{6}{5}$&$\frac{-18}{5}$&$\frac{1}{5}$&0&$\frac{6}{5}$&0&$\frac{-6}{5}$&$\frac{6}{5}$&$\frac{6}{5}$&$\frac{6}{5}$\\
&&&&&&&&&&&&\\
$4,6,1$&0&$\frac{-4}{5}$&0&$\frac{-6}{5}$&$\frac{24}{5}$&$\frac{-32}{5}$&$\frac{-36}{5}$&$\frac{48}{5}$&$\frac{24}{5}$&0&$\frac{-18}{5}$&$\frac{-24}{5}$\\
\hline
\end{tabular}
}
\end{center}

\newpage

\vglue 1.5cm 

\begin{center}
{\tiny
\textbf{Table 3 for the character} $\chi_{12}$.\\

\begin{tabular}{|c|c|c|c|c|c|c|c|c|c|c|c|c|c|c|}
\hline
&&&&&&&&&&&&&&\\
$a_1,a_2,b_1$ &$\!\!\beta_{1,\chi_{12}}\!\!$ & $\!\!\beta_{2,\chi_{12}}\!\!$ & $\!\!\beta_{3,\chi_{12}}\!\!$ & $\!\!\beta_{4,\chi_{12}}\!\!$ & $\!\!\beta_{5,\chi_{12}}\!\!$ & $\!\!\beta_{6,\chi_{12}}\!\!$ & $ \!\!\beta_{7,\chi_{12}}\!\!$ & $\!\!\beta_{8,\chi_{12}}\!\!$ & $\!\!\beta_{9,\chi_{12}}\!\!$ & $\!\!\beta_{10,\chi_{12}}\!\!$ & $\!\!\beta_{11,\chi_{12}}\!\!$ & $\!\!\beta_{12,\chi_{12}}\!\!$ & $\!\!\beta_{13,\chi_{12}}\!\!$ & $\!\!\beta_{14,\chi_{12}}\!\!$\\ 
&&&&&&&&&&&&&&\\
\hline 
&&&&&&&&&&&&&&\\       
$1,1,1$&-1&0&0&12&0&0&-4&0&0&3&0&0&0&0\\
&&&&&&&&&&&&&&\\  
$1,1,2$&-1&0&0&6&0&0&2&0&0&-3&0&0&0&0\\
&&&&&&&&&&&&&&\\  
$1,1,4$&-1&0&0&3&0&0&-1&0&0&3&0&0&0&0\\
&&&&&&&&&&&&&&\\  
$1,1,8$&$\frac{1}{2}$&$\frac{-3}{2}$&0&$\frac{3}{2}$&0&0&$\frac{1}{2}$&0&0&$\frac{3}{2}$&$\frac{9}{2}$&0&0&0\\
&&&&&&&&&&&&&&\\  
$1,1,16$&$\frac{-1}{4}$&$\frac{3}{4}$&$\frac{-3}{2}$&$\frac{3}{4}$&0&0&$\frac{-1}{4}$&0&0&$\frac{3}{4}$&$\frac{9}{4}$&$\frac{9}{2}$&3&0\\
&&&&&&&&&&&&&&\\  
$1,4,1$&$\frac{-1}{2}$&$\frac{1}{2}$&-1&6&-3&12&-2&-1&-4&$\frac{3}{2}$&$\frac{3}{2}$&3&3&-3\\
&&&&&&&&&&&&&&\\  
$1,4,2$&$\frac{-1}{2}$&$\frac{1}{2}$&-1&3&3&-12&1&-1&-4&$\frac{-3}{2}$&$\frac{-3}{2}$&-3&0&0\\
&&&&&&&&&&&&&&\\  
$1,4,4$&$\frac{-1}{2}$&$\frac{1}{2}$&-1&$\frac{3}{2}$&-3&12&$\frac{-1}{2}$&-1&-4&$\frac{3}{2}$&$\frac{3}{2}$&3&0&0\\
&&&&&&&&&&&&&&\\  
$1,4,8$&$\frac{1}{4}$&$\frac{-1}{4}$&-1&$\frac{3}{4}$&$\frac{-3}{2}$&6&$\frac{1}{4}$&$\frac{1}{2}$&2&$\frac{3}{4}$&$\frac{3}{4}$&-3&0&0\\
&&&&&&&&&&&&&&\\  
$1,4,16$&$\frac{-1}{8}$&$\frac{1}{8}$&-1&$\frac{3}{8}$&$\frac{-3}{4}$&3&$\frac{-1}{8}$&$\frac{-1}{4}$&-1&$\frac{3}{8}$&$\frac{3}{8}$&3&$\frac{3}{2}$&0\\
&&&&&&&&&&&&&&\\  
$2,2,1$&0&-1&0&9&-12&0&-3&-4&0&0&-3&0&0&0\\
&&&&&&&&&&&&&&\\  
$3,3,1$&-1&0&0&4&0&0&4&0&0&-1&0&0&0&0\\
&&&&&&&&&&&&&&\\  
$3,3,2$&-1&0&0&2&0&0&-2&0&0&1&0&0&0&0\\
&&&&&&&&&&&&&&\\  
$3,3,4$&-1&0&0&1&0&0&1&0&0&-1&0&0&0&0\\
&&&&&&&&&&&&&&\\  
$3,3,8$&$\frac{1}{2}$&$\frac{-3}{2}$&0&$\frac{1}{2}$&0&0&$\frac{-1}{2}$&0&0&$\frac{-1}{2}$&$\frac{-3}{2}$&0&0&0\\
&&&&&&&&&&&&&&\\  
$3,3,16$&$\frac{-1}{4}$&$\frac{3}{4}$&$\frac{-3}{2}$&$\frac{1}{4}$&0&0&$\frac{1}{4}$&0&0&$\frac{-1}{4}$&$\frac{-3}{4}$&$\frac{-3}{2}$&0&1\\
&&&&&&&&&&&&&&\\  
$3,12,1$&$\frac{-1}{2}$&$\frac{1}{2}$&-1&2&-1&4&2&1&4&$\frac{-1}{2}$&$\frac{-1}{2}$&-1&3&1\\
&&&&&&&&&&&&&&\\  
$3,12,2$&$\frac{-1}{2}$&$\frac{1}{2}$&-1&1&1&-4&-1&1&4&$\frac{1}{2}$&$\frac{1}{2}$&1&0&0\\
&&&&&&&&&&&&&&\\  
$3,12,4$&$\frac{-1}{2}$&$\frac{1}{2}$&-1&$\frac{1}{2}$&-1&4&$\frac{1}{2}$&1&4&$\frac{-1}{2}$&$\frac{-1}{2}$&-1&0&0\\
&&&&&&&&&&&&&&\\  
$3,12,8$&$\frac{1}{4}$&$\frac{-1}{4}$&-1&$\frac{1}{4}$&$\frac{-1}{2}$&2&$\frac{-1}{4}$&$\frac{-1}{2}$&-2&$\frac{-1}{4}$&$\frac{-1}{4}$&1&0&0\\
&&&&&&&&&&&&&&\\  
$3,12,16$&$\frac{-1}{8}$&$\frac{1}{8}$&-1&$\frac{1}{8}$&$\frac{-1}{4}$&1&$\frac{1}{8}$&$\frac{1}{4}$&1&$\frac{-1}{8}$&$\frac{-1}{8}$&-1&0&$\frac{1}{2}$\\
&&&&&&&&&&&&&&\\  
$4,4,1$&0&0&-1&$\frac{9}{2}$&-9&12&$\frac{-3}{2}$&-3&-4&0&0&3&3&-3\\
&&&&&&&&&&&&&&\\  
$6,6,1$&0&-1&0&3&-4&0&3&4&0&0&1&0&0&0\\
&&&&&&&&&&&&&&\\  
$12,12,1$&0&0&-1&$\frac{3}{2}$&-3&4&$\frac{3}{2}$&3&4&0&0&-1&3&1\\
&&&&&&&&&&&&&&\\ 
\hline
\end{tabular}
}
\end{center} 

\newpage

\begin{center}
{\tiny
\textbf{Table 3 for the character} $\chi_{24}$.\\

\begin{tabular}{|c|c|c|c|c|c|c|c|c|c|c|c|c|}
\hline
&&&&&&&&&&&&\\
$a_1,a_2,b_1$& $\!\!\beta_{1,\chi_{24}}\!\!$ & $\!\!\beta_{2,\chi_{24}}\!\!$ & $\!\!\beta_{3,\chi_{24}}\!\!$ & $\!\!\beta_{4,\chi_{24}}\!\!$ & $\!\!\beta_{5,\chi_{24}}\!\!$ & $\!\!\beta_{6,\chi_{24}}\!\!$ & $ \!\!\beta_{7,\chi_{24}}\!\!$ & $\!\!\beta_{8,\chi_{24}}\!\!$ & $\!\!\beta_{9,\chi_{24}}\!\!$ & $\!\!\beta_{10,\chi_{24}}\!\!$ & $\!\!\beta_{11,\chi_{24}}\!\!$ & $\!\!\beta_{12,\chi_{24}}\!\!$\\ 
&&&&&&&&&&&&\\
\hline 
&&&&&&&&&&&&\\  
$1,2,1$&$\frac{-1}{3}$&0&8&0&$\frac{8}{3}$&0&-1&0&0&0&$\frac{-4}{3}$&0\\
 &&&&&&&&&&&&\\  
$1,2,2$&$\frac{-1}{3}$&0&4&0&$\frac{-4}{3}$&0&1&0&-4&0&$\frac{-4}{3}$&0\\
&&&&&&&&&&&&\\  
$1,2,4$&$\frac{-1}{3}$&0&2&0&$\frac{2}{3}$&0&-1&0&0&0&$\frac{2}{3}$&0\\
&&&&&&&&&&&&\\  
$1,2,8$&$\frac{-1}{3}$&0&1&0&$\frac{-1}{3}$&0&1&0&2&0&$\frac{2}{3}$&0\\
&&&&&&&&&&&&\\  
$1,2,16$&$\frac{1}{6}$&$\frac{-1}{2}$&$\frac{1}{2}$&0&$\frac{1}{6}$&0&$\frac{1}{2}$&$\frac{3}{2}$&0&6&$\frac{2}{3}$&2\\
&&&&&&&&&&&&\\  
$2,4,1$&0&$\frac{-1}{3}$&6&-8&2&$\frac{8}{3}$&0&1&0&-16&-2&$\frac{-16}{3}$\\
&&&&&&&&&&&&\\  
$3,6,1$&$\frac{-1}{3}$&0&$\frac{8}{3}$&0&$\frac{8}{3}$&0&$\frac{-1}{3}$&0&$\frac{8}{3}$&0&$\frac{4}{3}$&0\\
&&&&&&&&&&&&\\  
$3,6,2$&$\frac{-1}{3}$&0&$\frac{4}{3}$&0&$\frac{-4}{3}$&0&$\frac{1}{3}$&0&$\frac{4}{3}$&0&0&0\\
&&&&&&&&&&&&\\  
$3,6,4$&$\frac{-1}{3}$&0&$\frac{2}{3}$&0&$\frac{2}{3}$&0&$\frac{-1}{3}$&0&$\frac{-4}{3}$&0&$\frac{-2}{3}$&0\\
&&&&&&&&&&&&\\  
$3,6,8$&$\frac{-1}{3}$&0&$\frac{1}{3}$&0&$\frac{-1}{3}$&0&$\frac{1}{3}$&0&$\frac{-2}{3}$&0&0&0\\
&&&&&&&&&&&&\\  
$3,6,16$&$\frac{1}{6}$&$\frac{-1}{2}$&$\frac{1}{6}$&0&$\frac{1}{6}$&0&$\frac{1}{6}$&$\frac{1}{2}$&$\frac{-4}{3}$&-2&$\frac{-2}{3}$&0\\
&&&&&&&&&&&&\\  
$6,12,1$&0&$\frac{-1}{3}$&2&$\frac{-8}{3}$&2&$\frac{8}{3}$&0&$\frac{1}{3}$&4&$\frac{16}{3}$&2&0\\
&&&&&&&&&&&&\\      

\hline
\end{tabular}
}
\end{center}

\section*{Acknowledgements}
We have used the open-source mathematics software SAGE (www.sagemath.org) for carrying out our calculations. 
Part of the work was done when the third named author was visiting NISER, Bhubaneswar and he thanks the institute for the 
warm hospitality.

\end{document}